\documentclass[12pt]{amsart}

\usepackage{epsfig}
\usepackage{amsmath,amsthm,amssymb,amsfonts}
\usepackage{mathrsfs,bbm}
\usepackage{xspace}
\usepackage[all,cmtip]{xy}

\newtheoremstyle{thm}{}{}{\itshape}{}{\scshape}{.}{0.14cm}{}
\theoremstyle{thm}
\newtheorem{theorem}{Theorem}[section]
\newtheorem{lemma}[theorem]{Lemma}
\newtheorem{proposition}[theorem]{Proposition}
\newtheoremstyle{def}{}{}{}{}{\scshape}{.}{0.2cm}{}
\theoremstyle{def}
\newtheorem{definition}[theorem]{Definition}
\newtheorem{example}[theorem]{Example}
\newenvironment{myproof}[1][Proof.]{\begin{trivlist}\item[\hskip \labelsep {\scshape #1}]}{\hfill\rule{1.05ex}{1.7ex}\end{trivlist}}

\newcommand{\cD}{\mathscr{D}}
\newcommand{\cF}{\mathscr{F}}
\newcommand{\cR}{\mathscr{R}}
\newcommand{\cS}{\mathscr{S}}
\newcommand{\PP}{P}
\newcommand{\E}{\mathrm{E}}  
\newcommand{\RR}{\mathbb{R}} 
\newcommand{\NN}{\mathbb{N}} 
\newcommand{\HH}{\mathbb{H}}
\newcommand{\I}{I}           
\newcommand{\SD}{S_\Delta}
\renewcommand{\phi}{\varphi}
\newcommand{\RN}{Radon-Nikodym\xspace}

\begin{document}

\title[Simple Random Densities]{Bayesian analysis of simple random densities}

\author[Marques]{Paulo C. Marques F.}
\address{Instituto de Matem\'atica e Estat\'istica da Universidade de S\~ao Paulo}
\email{pmarques@ime.usp.br}
\author[Pereira]{Carlos A. de B. Pereira}
\date{September 22, 2012}

\begin{abstract}
A tractable nonparametric prior over densities is introduced which is closed under sampling and exhibits proper posterior asymptotics.
\end{abstract}


\maketitle

\section{Introduction}
\thispagestyle{empty}

The early 1970's witnessed Bayesian inference going nonparametric with the introduction of statistical models with infinite dimensional parameter spaces; the most conspicuous being the Dirichlet process \cite{ferguson}, which is a prior on the class of all probability measures over a given sample space that trades great analytical tractability for a reduced support: as shown by Blackwell \cite{blackwell}, its realizations are, almost surely, discrete probability measures. The posterior expectation of a Dirichlet process is a probability measure that gives positive mass to each observed value in the sample, making the plain Dirichlet process unsuitable to handle inferential problems such as density estimation. Many extensions and alternatives to the Dirichlet process have been proposed \cite{ghosh}.

In this paper we construct a prior distribution over the class of densities with respect to Lebesgue measure. Given a partition in subintervals of a bounded interval of the real line, we define a random density whose realizations have a constant value on each subinterval of the partition. The distribution of the values of the random density on each subinterval is specified by transforming and conditioning a multivariate normal distribution.

Our simple random density is the finite dimensional analogue of the stochastic processes introduced by Thorburn \cite{thorburn} and Lenk \cite{lenk}.  Computations with these stochastic processes involve an intractable normalization constant, and are restricted to values of the random density on a finite number of arbitrarily chosen domain points, demanding some kind of interpolation of the results. The finite dimensionality of our random density makes our computations more direct and transparent and gives us simpler statements and proofs.

An outline of the paper is as follows. In Section \ref{sec:srd}, we give the formal definition of a simple random density. In Section \ref{sec:model}, we prove that the distribution of a simple random density is closed under sampling. The results of the simulations in Section \ref{sec:sims} depict the asymptotic behavior of the posterior distribution. We  extend the model hierarchically in Section \ref{sec:rp} to deal with random partitions. Although the usual Bayes estimate of a simple random density is a discontinuous density, in Section \ref{sec:smooth} we compute smooth estimates solving a decision problem in which the states of nature are realizations of the simple random density and the actions are smooth densities of a suitable class. Additional propositions and proofs of all the results in the paper are given in Section \ref{sec:proofs}.

\section{Simple random densities}\label{sec:srd}

Let $(\Omega,\cF,\PP)$ be the probability space from which we induce the distributions of all random objects considered in the paper. For some integer $k\geq 1$, let $\RR^k_+$ be the subset of vectors of $\RR^k$ with positive components. Write $\cR^k$ for the Borel sigma-field of $\RR^k$. Let $\lambda_k$ denote Lebesgue measure over $(\RR^k,\cR^k)$. We omit the indexes when $k=1$. The components of a vector $v\in\RR^k$ are written as $v_1,\dots,v_k$. 

Suppose that we have been given an interval $[a,b]\subset\RR$, and a set of real numbers $\Delta=\{t_0,t_1,\dots,t_k\}$, such that $a=t_0<t_1<\dots<t_k=b$, inducing a partition of $[a,b]$ into the $k\geq 1$ subintervals 
$$
  [a,t_1), [t_1,t_2), \dots, [t_{k-2},t_{k-1}), [t_{k-1},b] \, .
$$
The class of simple densities with respect to this partition consists of the nonnegative simple functions which have a constant value on each subinterval and integrate to one. Let $d_i=t_i-t_{i-1}$, for $i=1,\dots,k$, and define the map $\SD:\RR^k\to\RR$ by $\SD(u) = \sum_{i=1}^k d_i u_i$. Each simple density $f:\RR\to\RR$ within this class can be represented as
\[
  f(x)=\sum_{i=1}^k h_i\,\I_{[t_{i-1},t_i)}(x) \, ,
\]
in which $h=(h_1,\dots,h_k)\in\RR^k$ is such that each $h_i\geq 0$, and $\SD(h)=1$. The $h_i$'s are the heights of the steps of the simple density $f$.

From now on, let $\HH_r=\{v\in\RR^k_+:d_1 v_1+\dots+d_k v_k=r\}$, for $r\in\RR$. Note that, by the definition of the $d_i$'s given above, it follows that $\HH_r=\emptyset$ if $r\leq 0$. Moreover, define the projection on the first $k-1$ coordinates $\pi:\RR^k\to\RR^{k-1}$ by $\pi(v_1,\dots,v_{k-1},v_k)=(v_1,\dots,v_{k-1})$. For a normal random vector $Z=(Z_1,\dots,Z_k)$  with mean $m\in\RR^k$ and $k\times k$ covariance matrix $\Sigma$, denote by $U\sim L_k(m,\Sigma)$ the distribution of the lognormal random vector $U=(e^{Z_1},\dots,e^{Z_k})$. If $\Sigma$ is nonsingular, it is easy to show that $U$ has a density
\begin{align*}
  f_U(u) = (2\pi)^{-k/2} \, \vert\Sigma\vert^{-1/2} & \left(\prod_{i=1}^k u_i^{-1}\right) \\
  & \times \exp \left( -\frac{1}{2} (\log u-m)^\top\Sigma^{-1}\,(\log u-m)\right) \I_{\RR^k_+}(u) \, ,
\end{align*}
in which $\vert\Sigma\vert$ is the determinant of $\Sigma$, and we have introduced the notations $\log u = (\log u_1,\dots,\log u_k)^\top$ and $m=(m_1,\dots,m_k)^\top$. 

We define a random density whose realizations are simple densities with respect to the partition induced by $\Delta$ by specifying the distribution of the random vector of its steps heights. Informally, the steps heights will have the distribution of a lognormal random vector $U$ given that $\SD(U)=1$. The formal definition of the random density is given in terms of a version of the conditional distribution of $U$ given $\SD(U)$ and the expression of its conditional density with respect to a dominating measure. However, we are outside the elementary case in which the joint distribution is dominated by a product measure. In fact, we have in Proposition \ref{prop:singularity} a simple proof that Lebesgue measure $\lambda_{k+1}$ and the joint distribution of $U$ and $\SD(U)$ are mutually singular. 

A suitable family of measures that dominate the conditional distribution of $U$ given $\SD(U)$, for each value of $\SD(U)$, is described in the following lemma.

\begin{lemma}\label{lemma:tau}
Let $\tau_r:\cR^k\to\RR$ be defined by $\tau_r(A)=d_k^{-1} \lambda_{k-1}(\pi(A\cap\HH_r))$, for $r\in\RR$. Then, each $\tau_r$ is a measure over $(\RR^k,\cR^k)$. 
\end{lemma}

\begin{figure}[t!] 
\begin{center}
\includegraphics[width=8.5cm]{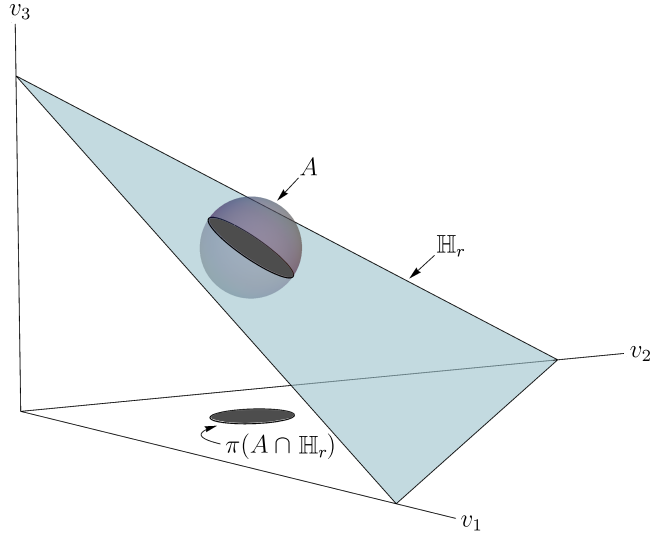}
\end{center}
\caption{Geometrical interpretation of the measures $\tau_r$ of Lemma \ref{lemma:tau}, for $r>0$, in the particular case when $k=3$. The value of $\tau_r(A)$ is the area of the projection $\pi(A\cap\HH_r)$ multiplied by $d_3^{-1}$.}
\label{fig:geo} 
\end{figure}

The proof of Lemma \ref{lemma:tau} is given in Section \ref{sec:proofs}. Figure \ref{fig:geo} gives a simple geometric interpretation of the measures $\tau_r$ when the underlying partition is formed by three subintervals.

The following result is the basis for the formal definition of the random density.

\begin{theorem}\label{theo:cond.dens}
Let $U\sim L_k(m,\Sigma)$, with nonsingular $\Sigma$, and let $\{\tau_r\}_{r\in\RR}$ be the family of measures over $(\RR^k,\cR^k)$ defined on Lemma \ref{lemma:tau}. Then, we have that $\mu_{U\mid\SD(U)}:\cR^k\times\RR_+\to\RR$ defined by
\[
  \mu_{U\mid\SD(U)}(A\mid r) = \int_A \frac{f_U(u)}{f_{\SD(U)}(r)}\,\I_{\HH_r}(u)\,d\tau_r(u) \, ,
\]
is a regular version of the conditional distribution of $U$ given $\SD(U)$, in which
\[
  f_{\SD(U)}(r) = {\displaystyle \int_{\RR^k}} f_U(u)\,\I_{\HH_r}(u)\,d\tau_r(u) \, .
\]
Moreover, $\mu_{U\mid\SD(U)}(\HH_r\mid r)=1$, for each $r>0$.
\end{theorem}

The necessary lemmata and the proof of Theorem \ref{theo:cond.dens} are given in Section \ref{sec:proofs}. The following definition of the random density uses the specific version of the conditional distribution constructed in Theorem \ref{theo:cond.dens}.

\begin{definition}\label{def:srd}
Let $U\sim L_k(m,\Sigma)$, with nonsingular $\Sigma$. We say that the map $\phi:\RR\times\Omega\to\RR$ defined by
\[
  \phi(x,\omega) = \sum_{i=1}^k H_i(\omega)\,\I_{[t_{i-1},t_i)}(x)
\] 
is a \textit{simple random density}, in which $H=(H_1,\dots,H_k)$ are the \textit{random heights of the steps} of $\phi$, with distribution given by $\mu_H(A)=\mu_{U\mid \SD(U)}(A\mid 1)$, for $A\in\cR^k$, and $\mu_{U\mid \SD(U)}$ is the regular version of the conditional distribution of $U$ given $\SD(U)$ obtained in Theorem \ref{theo:cond.dens}. Hence, for every $A\in\cR^k$, we have
\[
  \mu_H(A) = \int_A \frac{f_U(h)}{f_{\SD(U)}(1)}\,\I_{\HH_1}(h)\,d\tau_1(h) \, ,
\]
in which $\tau_1(A)=d_k^{-1} \lambda_{k-1}(\pi(A\cap\HH_1))$ and it holds that $\mu_H(\HH_1)=1$. We use the notation $\phi\sim\Delta(m,\Sigma)$.
\end{definition}

\section{Conditional model}\label{sec:model}

Now, we model a set of absolutely continuous observables conditionally, given the value of a simple random density $\phi$. The following lemma, proved in Section \ref{sec:proofs}, describes the conditional model and determines the form of the likelihood.

\begin{lemma}\label{lemma:model}
Consider $\phi\sim\Delta(m,\Sigma)$ represented as
$$
  \phi(x,\omega) = \sum_{i=1}^k H_i(\omega)\,\I_{[t_{i-1},t_i)}(x) \, ,
$$ 
and let the random variables $X_1,\dots,X_n$ be conditionally independent and identically distributed, given that $H=h$, with distribution 
$$
  \mu_{X_1\mid H}(A\mid h)=\int_A f(y)\,d\lambda(y)\, ,
$$
in which we have defined $f(y)=\sum_{i=1}^k h_i\,\I_{[t_{i-1},t_i)}(y)$. Define $X=(X_1,\dots,X_n)$ and let $x=(x_1,\dots,x_n)\in\RR^n$. Then, $\mu_{X\mid H}(\;\cdot\mid h)\ll\lambda_n$, almost surely $[\mu_H]$, with \RN derivative
\[
  \frac{d\mu_{X\mid H}}{d\lambda_n}(x\mid h) = f_{X\mid H}(x\mid h) = \prod_{i=1}^k h_i^{c_i} \, ,
\]
in which $c_i=\sum_{j=1}^n \I_{[t_{i-1},t_i)}(x_j)$, for $i=1,\dots,k$.
\end{lemma}

The factorization criterion implies that $c=(c_1,\dots,c_n)$ is a sufficient statistic for $\phi$. That is, in this conditional model, as one should expect, all the sample information is contained in the countings of how many sample points belong to each subinterval of the partition induced by $\Delta$.

Using the notation of Lemma \ref{lemma:model}, and defining $c=(c_1,\dots,c_k)^\top$, we can prove that the prior distribution of $\phi$ is closed under sampling.

\begin{theorem}\label{theo:closure}
If $\phi\sim\Delta(m,\Sigma)$, then $\phi\mid X=x\sim\Delta(m^*,\Sigma)$, in which $m^*=m+\Sigma c$.
\end{theorem}

This result, proved in Section \ref{sec:proofs}, makes the simulations of the prior and posterior distributions essentially the same, the only difference being the computation of $m^*$.

\section{Stochastic simulations}\label{sec:sims}

We summarize the distribution of a simple random density $\phi\sim\Delta(m,\Sigma)$, represented as $\phi(x,\omega) = \sum_{i=1}^k H_i(\omega)\,\I_{[t_{i-1},t_i)}(x)$, in two ways. First, motivated by the fact, proved in Proposition \ref{prop:pred}, that the prior and posterior expectations are predictive densities, we take as an estimate the expectation of the steps heights $\hat{h}=(\E[H_1],\dots,\E[H_k])$. Second, the uncertainty of this estimate is assessed defining
\[
  B(\hat{h},\epsilon) = \left\{ h\in\HH_1 : \max_{1\leq i\leq k}|\hat{h}_i - h_i| < \epsilon \right\} \, ,
\]
for $\epsilon>0$, and taking as a credible set the $B(\hat{h},\epsilon)$ with the smallest $\epsilon$ such that $\PP\{\omega:H(\omega)\in B(\hat{h},\epsilon)\}=\gamma$, in which $\gamma\in(0,1)$ is the credibility level.

The Random Walk Metropolis algorithm \cite{robert} is used to draw dependent realizations of the steps of $\phi$ as values of a Markov chain $\{H^{(i)}\}_{i\geq 0}$. The two summaries are computed through ergodic means of this chain. For example, the credible set is determined with the help of the almost sure convergence of
\[
  \frac{1}{N} \sum_{i=0}^N \I_{B(\hat{h},\epsilon)} \left(H^{(i)}\right) \xrightarrow[N\to\infty]{} \E\left[\I_{B(\hat{h},\epsilon)}(H)\right] = 	\PP\left\{\omega:H(\omega)\in B(\hat{h},\epsilon)\right\} \, .
\]

As for the parameters appearing in Definition \ref{def:srd}, we take in our experiments all the $m_i$'s equal to one, and the covariance matrix $\Sigma=(\sigma_{ij})$ is chosen in the following way. Given some positive definite covariance function $C:\RR\times\RR\to\RR$, we induce $\Sigma$ from $C$ defining
\[
  \sigma_{ij} = C\!\left( \frac{t_{i-1}+t_i}{2}, \frac{t_{j-1}+t_j}{2} \right)\, ,
\]
for $i,j=1,\dots,k$. In our examples we study the family of Gaussian covariance functions defined by $C_{\rho,\theta}(x,y)=\rho\,e^{-\theta\,(x-y)^2}$, with dispersion parameter $\rho>0$ and scale parameter $\theta>0$.

\begin{figure}[t!] 
\begin{center}
\includegraphics[width=12.5cm]{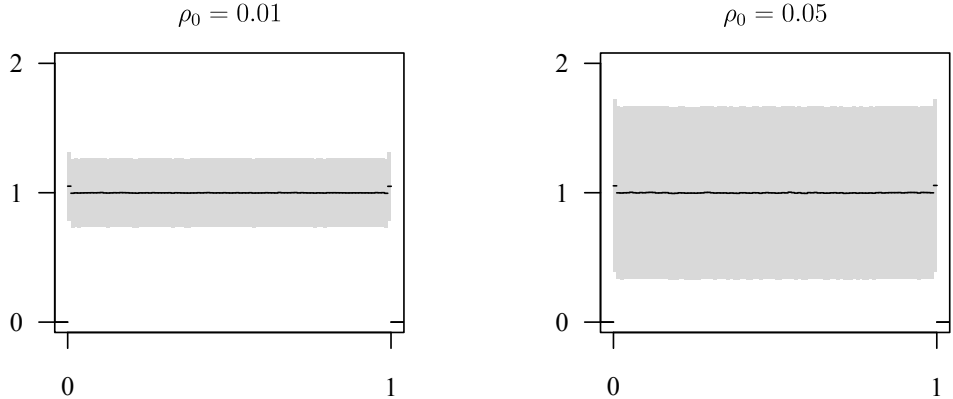}
\end{center}
\caption{Effect of the value of $\rho_0$ on the concentration of the prior. The curves in black are prior expectations and the gray regions are credible sets with credibility level of $95\%$.}
\label{fig:prior} 
\end{figure}

\begin{figure}[t!] 
\begin{center}
\includegraphics[width=12.5cm]{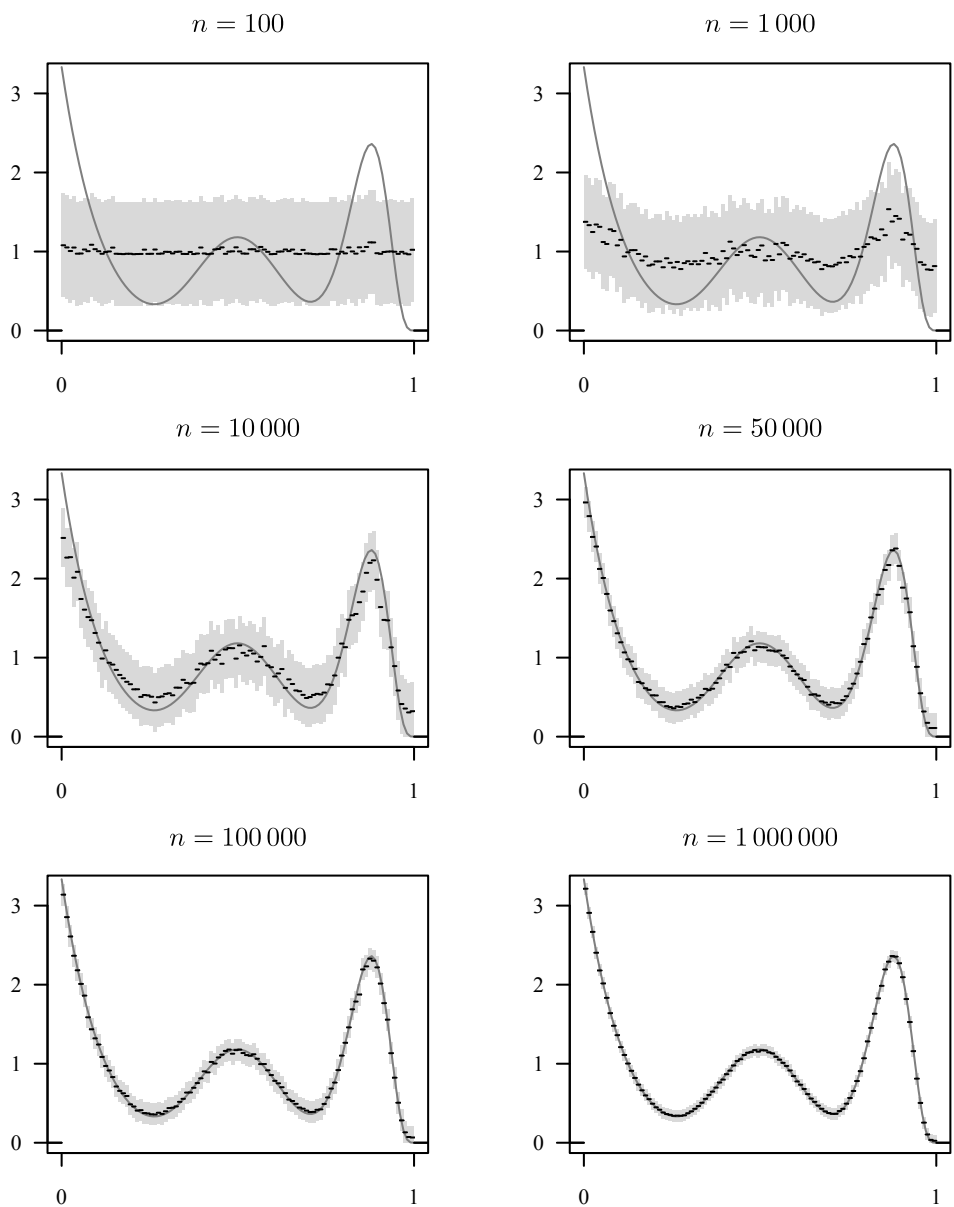}
\end{center}
\caption{Posterior summaries for Example \ref{ex:betas}. On each graph, the black simple density is the estimate $\hat{\phi}$, the light gray region is a credible set with credibility level of $95\%$, and the dark gray curve is the data generating density.}
\label{fig:betas} 
\end{figure}

\begin{example}\label{ex:betas}
Let $\phi\sim\Delta(m,\Sigma)$ and consider the sample space $[0,1]$ with $\Delta = \left\{ 0, 0.01, 0.02, \dots, 0.98, 0.99, 1 \right\}$. For the sake of generality, we induce $\Sigma$ from the family of Gaussian covariance functions with fixed dispersion parameter $\rho_0$ but with random scale parameter $\Theta=Y+20\,000$, in which $Y\sim \mathrm{Gamma}(2,0.001)$. These choices guarantee that computations with $\Sigma$ are numerically stable. In Figure \ref{fig:prior}, the summaries of the prior distribution of $\phi$ show that the value of $\rho_0$ controls the concentration of the prior. Fixing $\rho_0=0.05$ and generating data from the mixture
\[
  \frac{1}{3}\cdot \mathrm{Beta}(1,10) + \frac{1}{3}\cdot \mathrm{Beta}(10,10) + \frac{1}{3}\cdot \mathrm{Beta}(30,5) \, ,
\]
we have in Figure \ref{fig:betas} the posterior summaries for different sample sizes. Note the concentration of the posterior as we increase the size of the samples. \hfill\rule{1.05ex}{1.7ex}
\end{example}

We observe the same asymptotic behavior of the posterior distribution with data coming from a triangular distribution and a mixture of normals (with appropriate truncation of the sample space).

\section{Random partitions}\label{sec:rp}

Inferentially, we have a richer construction when the definition of the simple random density involves a random partition. Informally, we want a model for the random density in which the underlying partition adapts itself according to the information contained in the data.

We consider a family of uniform partitions of a given interval $[a,b]$. Each partition of this family will be described by a positive integer random variable $K$, which determines the number of subintervals in the partition. Since the parameter $\rho$ of the family of Gaussian covariance functions used to induce $\Sigma$ may have different meanings for different partitions, we treat it as a positive random variable $R$.

Explicitly, we are considering the following hierarchical model: $K$ and $R$ are independent. Given that $K=k$ e $R=\rho$, we choose the uniform partition of the interval $[a,b]$ induced by
\[
  \Delta = \left\{ a, a + \frac{b - a}{k}, a + \frac{2(b - a)}{k}, \dots, a + \frac{(k-1)(b - a)}{k}, b \right\} \, ,
\]
induce $\Sigma_{\rho,\theta}$ from the family of Gaussian covariance functions, and take the simple random density $\phi\sim \Delta(m,\Sigma_{\rho,\theta})$. Finally, the observables are modeled as in Lemma \ref{lemma:model}. This hierarchy is summarized in the following graph.
\[
  \entrymodifiers={++[o][F-]}
  \SelectTips{cm}{}
  \vcenter{\xymatrixcolsep{0.8cm}\xymatrixrowsep{0.5cm}\xymatrix{
    K \ar[dr] & *{} & *{} & R \ar[dll] \\
    *{} & \phi \ar[dl] \ar[d] \ar[drr] \\
    X_1 & X_2 & *{\dots} & X_n
  }}
\]

In the following example, instead of specifying priors for $K$ and $R$, we define the likelihood of $K$ and $R$ by $L_x(k,\rho) = f_{X\mid K,R}(x\mid k,\rho)$, whose form is obtained in Proposition \ref{prop:Lkr}, find the maximum $(\hat{k},\hat{\rho}) = \arg \max_{k,\rho} L_x(k,\rho)$, and use these values in the definitions of the prior, determining the posterior summaries as we did in Section \ref{sec:sims}. 

\begin{figure}[t!] 
\begin{center}
\includegraphics[width=10cm]{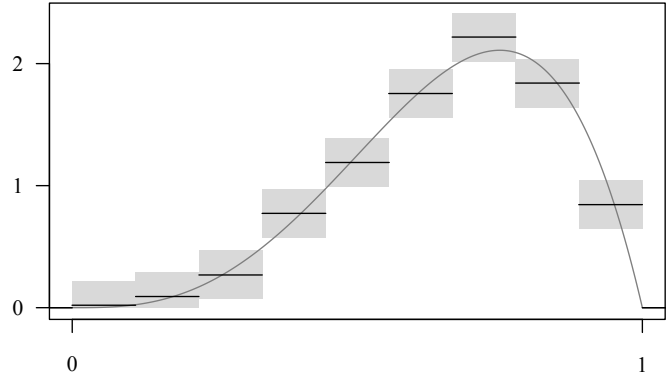}
\end{center}
\caption{Posterior summaries for Example \ref{ex:eb}. The black simple density is the estimate $\hat{\phi}$, the light gray region is a credible set with credibility $95\%$, and the dark gray curve is the data generating density.}
\label{fig:eb} 
\end{figure}

\begin{figure}[t!]
\begin{center}
  \begin{minipage}[c]{6.5cm}
    \includegraphics[width=5.8cm]{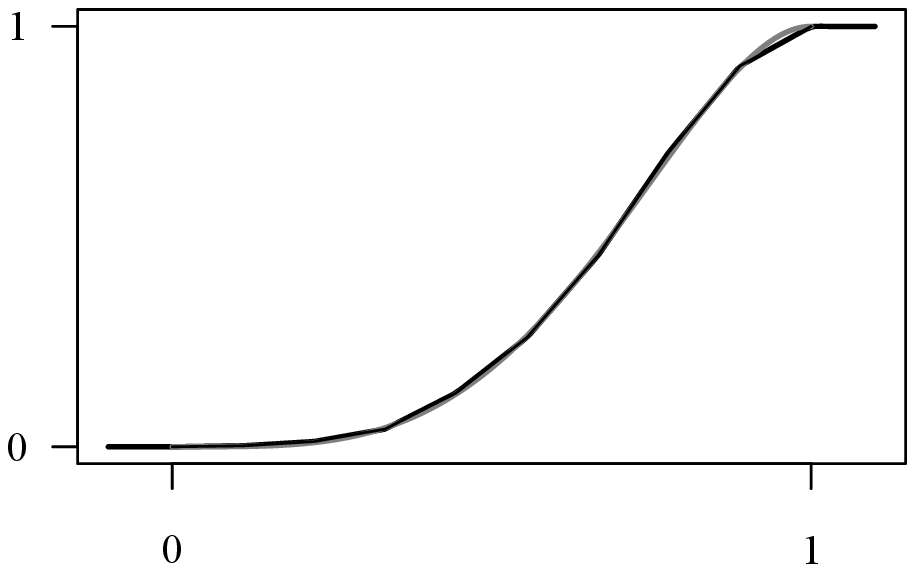}
  \end{minipage}
  \begin{minipage}[c]{6.5cm}
    \includegraphics[width=5.8cm]{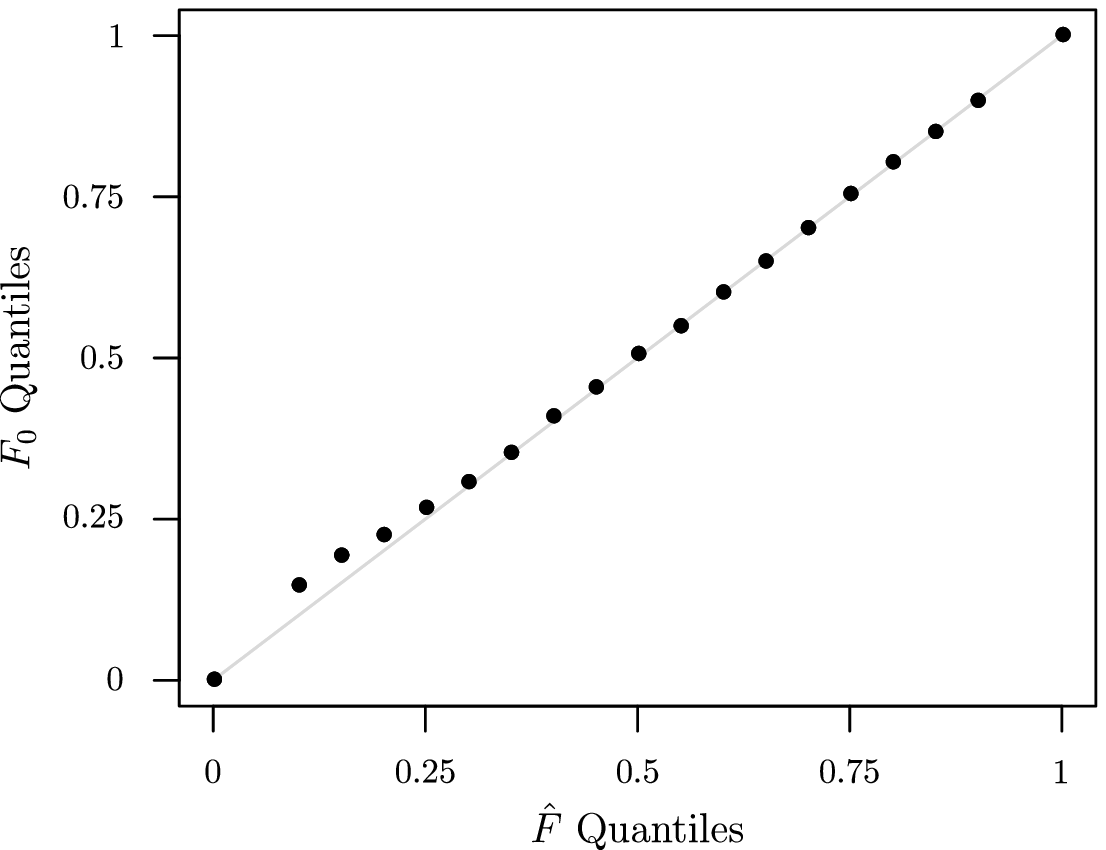}
  \end{minipage}
\end{center}
\caption{Example \ref{ex:eb}. On the left graph, the black curve is the estimated distribution function $\hat{F}$ and the gray curve is the data generating distribution function $F_0$. On the right graph, we have the comparison of some of the quantiles of $\hat{F}$ and $F_0$.} 
\label{fig:F.QQ}
\end{figure}

\begin{example}\label{ex:eb} 
With a sample of size $2\,000$ generated from a $\mathrm{Beta}(4,2)$ distribution, we find the maximum of the likelihood of $K$ and $R$ at $(\hat{k},\hat{\rho})=(9,1.43)$. In Figure \ref{fig:eb} we have the posterior summaries obtained using these values in the definition of the prior. Moreover, in the left graph of Figure \ref{fig:F.QQ} we have the distribution function $\hat{F}$ corresponding to the estimated posterior density. For the sake of comparison, we plot in the right graph of Figure \ref{fig:F.QQ} some quantiles of this distribution $\hat{F}$ against the quantiles of the distribution $F_0$ from which we generated the data. \hfill\rule{1.05ex}{1.7ex}
\end{example}

\section{Smooth estimates}\label{sec:smooth}

It is possible to go beyond the discontinuous densities obtained as estimates in the last two sections and get smooth estimates of a simple random density $\phi$ solving a Bayesian decision problem in which the states of nature are the realizations of $\phi$ and the actions are smooth densities of a suitable class. 

In view of Theorem \ref{theo:closure}, it is enough to consider the problem without data. As before, the sample space is the interval $[a,b]$, which is partitioned according to some $\Delta$. For a density $f$ with respect to Lebesgue measure, we denote its $L_2$ norm by $\Vert f\Vert_2 = \left( \int f^2 d\lambda \right)^{1/2}$.

\begin{proposition}\label{prop:decision}
For $N\geq 1$, let $g_1,\dots,g_N$ be densities with respect to Lebesgue measure, with support $[a,b]$, such that $\Vert g_i\Vert_2<\infty$, and let $\cD$ be the class of densities of the form $\sum_{i=1}^N \alpha_i\,g_i$, with $\alpha_i\geq 0$, for $i=1,\dots,N$, and $\sum_{i=1}^N \alpha_i=1$. Let $\phi\sim\Delta(m,\Sigma)$ and define $\cS$ as the class of densities which are realizations of $\phi$. Define the loss function $L:\cS\times\cD\to\RR$ by
\[
  L(s,f)=\Vert s-f\Vert_2^2=\int_a^b \left(s(x)-f(x)\right)^2\,d\lambda(x) \, .
\]
Then, the Bayes decision is $\hat{\phi}=\sum_{i=1}^N \hat{\alpha}_i\,g_i$, in which $\hat{\alpha_i}$ minimize globally the quadratic form
\[
  Q = \sum_{i,j=1}^N \alpha_i \alpha_j \, M_{ij} - \sum_{i=1}^N \alpha_i \, J_i \, ,
\]
subject to the constraints $\alpha_i\geq 0$, for $i=1,\dots,N$, and $\sum_{i=1}^N \alpha_i=1$, with the definitions
\[
    M_{ij} = \int_a^b g_i(x) g_j(x)\,d\lambda(x) \, ,\qquad J_i = 2 \int_a^b g_i(x) \E[\phi(x)]\,d\lambda(x) \, .
\]
\end{proposition}

We use the result of Proposition \ref{prop:decision}, proved in Section \ref{sec:proofs}, choosing the $g_i$'s inside a class of smooth densities that serve approximately as a basis to represent any continuous density with the specified support.

\begin{figure}[t!] 
\begin{center}
\includegraphics[width=12cm]{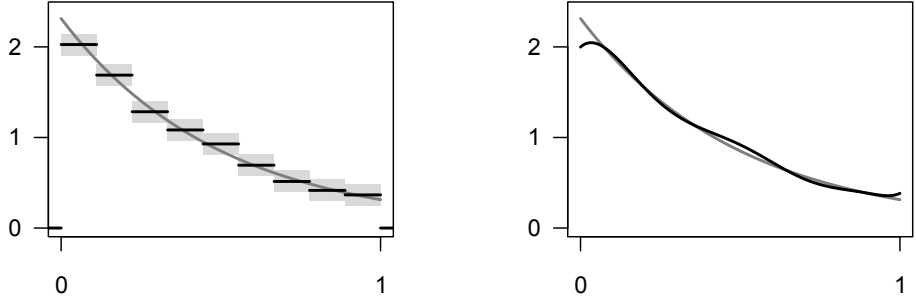}
\end{center}
\caption{Example \ref{ex:exp}. On the right graph, the black simple density is the estimate $\hat{\phi}$, and the light gray region is a credible set with credibility $95\%$. On both graphs the dark gray curve is the data generating density. On the left graph, the black smooth density is the Bayes decision of Proposition \ref{prop:decision}.}
\label{fig:smooth} 
\end{figure}

For the next example, suppose that the support of the densities is the interval $[0,1]$. Bernstein's Theorem (see \cite{billi}, Theorem 6.2) states that the polynomial 
\[
  B_N(x) = \sum_{i=0}^N f\!\left(\frac{i}{N}\right) \binom{N}{i} \, x^i (1-x)^{N-i}
\]
approximates uniformly any continuous function $f$ defined on $[0,1]$, when $N\to\infty$. Suppose that $f$ is a density. If we define, for $i=0,\dots,N$,
\[
  \alpha_i = f\!\left(\frac{i}{N}\right) \binom{N}{i} \frac{\Gamma(i+1)\Gamma(N-i+1)}{\Gamma(N+2)} \, ,
\]
we can rewrite the approximating polynomial as $B_N(x) = \sum_{i=0}^N \alpha_i\,g_i(x)$, in which $g_i$ is a density of a $\mathrm{Beta}(i+1,N-i+1)$ random variable. Hence, if we take a sufficiently large $N$, we expect that any continuous density with support $[0,1]$ will be reasonably approximated by a mixture of these $g_i$'s.

\begin{example}\label{ex:exp}
Suppose that we have a sample of $5\,000$ data points simulated from a truncated exponential distribution, whose density is 
\[
  f_0(x) = \frac{2\,e^{-2(x-1)}}{e^2-1} \,\,\I_{[0,1]}(x) \, .
\]
Repeating the analysis made in Example \ref{ex:eb}, we find the maximum of the likelihood of $K$ e $R$ at $(\hat{k},\hat{\rho}) = (9,0.86)$. The left graph of Figure \ref{fig:smooth} presents the posterior summaries. After that, we solved the problem of constrained optimization in Proposition \ref{prop:decision} and found the results shown in the right graph of Figure \ref{fig:smooth}. \hfill\rule{1.05ex}{1.7ex}
\end{example}

\section{Additional results and proofs}\label{sec:proofs}

In this section we present some subsidiary propositions and give proofs to all the results stated in the paper.

\begin{proposition}\label{prop:singularity}
Let $U\sim L_k(m,\Sigma)$ and denote by $\mu_{U,\SD(U)}$ the joint distribution of $U$ and $\SD(U)$. Then, $\mu_{U,\SD(U)}\perp\lambda_{k+1}$.
\end{proposition}

\begin{myproof}
Define the set $A = \left\{ v\in\RR^{k+1}:\sum_{i=1}^k d_i v_i = v_{k+1} \right\}\in\cR^{k+1}$. Then,
\begin{align*}
  \mu_{U,\SD(U)}(A) &= \PP\left\{\omega:(U(\omega),\SD(U(\omega)))\in A\right\} \\ 
  &= \PP\left\{\omega: \sum_{i=1}^k d_i U_i(\omega) = \SD(U(\omega)) \right\} = 1 ,
\end{align*}
by definition of $\SD$. On the other hand, note that $\lambda_{k+1}(A)=0$, since this is the $(k+1)$-volume of the $k$-dimensional hyperplane defined by the set $A$. Since $\mu_{U,\SD(U)}(A^c) = 0$, the result follows.
\end{myproof}

\begin{myproof}[Proof of Lemma \ref{lemma:tau}.]
When $r\leq 0$, the result is trivial, since in this case $\HH_r=\emptyset$, making $\tau_r$ a null measure. Suppose that $r>0$ and let $g:\RR^k\to\RR^k$ be the function defined by
\[
  g(v)=\left(v_1,\dots,v_{k-1},\frac{1}{d_k}\left(v_k-\sum_{i=1}^{k-1}d_i v_i\right)\right) \, .
\]
Define $h_r:\RR^{k-1}\to\RR^k$ by $h_r(y)=g(y,r)$. We will show that $\pi(A\cap\HH_r)=h^{-1}_r(A)$, for every $A\in\cR$. Suppose that $y\in\pi(A\cap\HH_r)$. Then, there is a $v\in A\cap\HH_r$ such that $y=\pi(v)=(v_1,\dots,v_{k-1})$ and
\[
  h_r(y) = g(y,r) = \left( v_1,\dots,v_{k-1}, \frac{1}{d_k}\left(r - \sum_{i=1}^{k-1}d_i v_i\right)\right) \, .
\]
Since $v\in\HH_r$, we have that $\frac{1}{d_k}\left(r-\sum_{i=1}^{k-1}d_i v_i\right)=v_k$, implying that $h_r(y)=v$. Since $v\in A$, it follows from the definition of the inverse image of $h_r$ that $y\in h_r^{-1}(A)$ and, therefore, we conclude that $\pi(A\cap\HH_r)\subset h_r^{-1}(A)$. To prove the other inclusion, suppose that $y\in h_r^{-1}(A)$ and define $v=h_r(y)$. Hence, $v\in A$ and by the definition of $h_r$ we have that
\[
  v = g(y,r) = \left( y_1,\dots,y_{k-1}, \frac{1}{d_k}\left(r - \sum_{i=1}^{k-1}d_i y_i\right)\right) \, ,
\]
implying that $v\in \HH_r$, because $\sum_{i=1}^k d_i v_i=r$. Since $v\in A\cap\HH_r$ and $y=\pi(v)$, it follows that $y\in\pi(A\cap\HH_r)$. Therefore, $h_r^{-1}(A)\subset \pi(A\cap\HH_r)$. Hence, we have that $\tau_r=d_k^{-1} \lambda_k\circ h_r^{-1}$ and the properties of the inverse image of $h_r$ and the Lebesgue measure entail that each $\tau_r$ is a measure over $(\RR^k,\cR^k)$.
\end{myproof}

\begin{lemma}\label{lemma:joint.dens}
Let $U\sim L_k(m,\Sigma)$. Let $\xi$, defined by 
$$
  \xi(A)=\lambda_k\{u\in\RR^k_+:(u,\SD(u))\in A\} \, ,
$$
be a measure over $(\RR^{k+1},\cR^{k+1})$. Denote by $\mu_{U,\SD(U)}$ the joint distribution of $U$ and $\SD(U)$. Then, we have that $\mu_{U,\SD(U)}\ll\xi$, with \RN derivative $d\mu_{U,\SD(U)}/d\xi=f_{U,\SD(U)}$ given by
\[
  f_{U,\SD(U)}(u,r) = f_U(u)\,\I_{\HH_r}(u) \, ,
\]
in which $u\in\RR^k$ and $r\in\RR$. 
\end{lemma}

\begin{myproof}
Define the function $T:\RR^k_+\to\RR^{k+1}$ by $T(u)=(u,\SD(u))$. Note that $\xi=\lambda_k\circ T^{-1}$. Define the function $\psi:\RR^{k+1}\to\RR$ by $\psi(u,r)=f_U(u)\,\I_{\HH_r}(u)$, with $u\in\RR^k$ and $r\in\RR$. The diagram
\[
  \xymatrixcolsep{1.7cm}\xymatrixrowsep{1.7cm}\xymatrix{\RR^k_+ \ar[rd]_{f_U} \ar[r]^T
            &\RR^{k+1} \ar[d]^\psi \\
            &\RR}
\]
commutes, since $\psi(T(u)) = \psi(u, \SD(u)) = f_U(u)\,\I_{\HH_{\SD(U)}}(u) = f_U(u)$, for every $u\in\RR^k_+$. For every $A\in\cR^{k+1}$, we have that
\begin{align*}
  \mu_{U,\SD(U)}(A) &= \PP\{\omega:(U(\omega),\SD(U(\omega)))\in A\} = \PP\{\omega:U(\omega)\in T^{-1}(A)\} \\
    &= \int_{T^{-1}(A)} f_U(u)\,d\lambda_k(u) = \int_{T^{-1}(A)} \psi(T(u))\,d\lambda_k(u) \\
    &= \int_A \psi(u,r)\,d\xi(u,r) = \int_A f_U(u)\,\I_{\HH_r}(u)\,d\xi(u,r) \, ,
\end{align*}
in which the fifth equality is obtained transforming by $T$, $u\in\RR^k$ and $r\in\RR$. It follows that $\mu_{U,\SD(U)}\ll\xi$, and the \RN derivative has the desired expression.
\end{myproof}

\begin{lemma}\label{lemma:desint}
Let $\xi$ be the measure defined on Lemma \ref{lemma:joint.dens} and let $\{\tau_r\}_{r\in\RR}$ be the family of measures defined on Lemma \ref{lemma:tau}. Then, for every measurable nonnegative $\psi:\RR^{k+1}\to\RR$, we have that
\[
  \int_{\RR^{k+1}} \psi(u,r)\,d\xi(u,r) = \int_{\RR}\left(\int_{\RR^k} \psi(u,r)\,d\tau_r(u)\right)d\lambda(r) \, ,
\]
in which $u\in\RR^k$ and $r\in\RR$.
\end{lemma}

\begin{myproof}
Define $f:\RR^k\to\RR^k$ by $f(u)=(u_1,\dots,u_{k-1},\sum_{i=1}^k d_i u_i)$. Hence, $f$ is a differentiable function whose inverse is the differentiable function $g$ defined on Lemma \ref{lemma:tau}. The value of the Jacobian on the point $v\in\RR^k$ is $J_g(v)=d_k^{-1}$. Let $A\in\cR^k$, $y\in\RR^{k-1}$, $r\in\RR$, and define $h_r$ as in Lemma \ref{lemma:tau}. When $r>0$, we have already shown in the course of the proof of Lemma $\ref{lemma:tau}$ that $\pi(A\cap\HH_r)=h_r^{-1}(A)$, for every $A\in\cR^k$. Remembering that, by definition, $\HH_r\subset\RR^k_+$, it follows that $\pi(A\cap\HH_r)=h_r^{-1}(A\cap\RR^k_+)$ and we conclude that $\I_{\pi(A\cap\HH_r)}(y)=\I_{A\cap\RR^k_+}(g(y,r))$. Now suppose that $r\leq 0$. In this case, since $\HH_r=\emptyset$, we have that $\I_{\pi(A\cap\HH_r)}(y)=\I_\emptyset(y)=0$. As for the value of $\I_{A\cap\RR^k_+}(g(y,r))$, consider two subcases: since
\[
  g(y,r) = \left( y_1,\dots,y_{k-1}, \frac{1}{d_k}\left(r - \sum_{i=1}^{k-1}d_i y_i\right)\right) \, ,
\]
if any of the $y_i\leq 0$, then $\I_{A\cap\RR^k_+}(g(y,r))=0$, otherwise, we have 
$$
  \frac{1}{d_k}\left(r - \sum_{i=1}^{k-1}d_i y_i\right)<0 \, ,
$$
and again it happens that $\I_{A\cap\RR^k_+}(g(y,r))=0$. Therefore, we conclude that in this case also $\I_{\pi(A\cap\HH_r)}(y)=\I_{A\cap\RR^k_+}(g(y,r))$. Hence, for $A\in\cR^k$ and $B\in\cR$, we have that
\begin{align*}
  \xi(A\times B) &= \lambda_k\{u\in\RR^k_+ : u\in A, \SD(u)\in B\} \\
                 &= \int_{\RR^k} \I_{A\cap\RR^k_+}(u)\,\I_B(\SD(u))\,d\lambda_k(u) \\
                 &= \int_{\RR^k} \I_{A\cap\RR^k_+}(g(y,r))\,\I_B(r)\,\vert J_g(y,r)\vert\,d\lambda_k(y,r) \\
                 &= \int_{\RR^k} d_k^{-1} \, \I_{\pi(A\cap\HH_r)}(y) \I_B(r)\,d\lambda_k(y,r) \\
                 &= \int_B\left(d_k^{-1} \int_{\pi(A\cap\HH_r)} d\lambda_{k-1}(y)\right)d\lambda(r) = \int_B \tau_r(A)\,d\lambda(r) \, ,
\end{align*}
in which $y\in\RR^{k-1}$ and $r\in\RR$, the third equality is obtained transforming by $f$, and the penultimate equality is a consequence of Tonelli's Theorem. The result follows from the Product Measure Theorem and Fubini's Theorem (see \cite{ash}, Theorems 2.6.2 e 2.6.4).
\end{myproof}

\begin{lemma}\label{lemma.dens.S}
Let $U\sim L_k(m,\Sigma)$. Let $\{\tau_r\}_ {r\in\RR}$ be the family of measures defined on Lemma \ref{lemma:tau}. Let $\mu_{\SD(U)}$ be the distribution of $\SD(U)$. Then, $\mu_{\SD(U)}\ll\lambda$ with \RN derivative $d\mu_{\SD(U)}/d\lambda=f_{\SD(U)}$ given by 
$$
  f_{\SD(U)}(r) = {\displaystyle\int_{\RR^k}} f_U(u)\,\I_{\HH_r}(u)\,d\tau_r(u) \, .
$$
\end{lemma}

\begin{myproof}
Let $A\in\cR$, $u\in\RR^k$, and $r\in\RR$. Let $\xi$ be the measure defined on Lemma \ref{lemma:joint.dens}. We have that
\begin{align*}
  \mu_{\SD(U)}(A) &= \PP\{\omega : \SD(U(\omega))\in A\} = \PP\{\omega : U(\omega)\in\RR^k, \SD(U(\omega))\in A\} \\
                  &= \mu_{U,\SD(U)}(\RR^k\times A) = \int_{\RR^k\times A} f_U(u)\,\I_{\HH_r}(u)\,d\xi(u,r) \\
                  &= \int_A\left(\int_{\RR^k} f_U(u)\,\I_{\HH_r}(u)\,d\tau_r(u)\right) d\lambda(r) \, ,
\end{align*}
in which the penultimate equality follows from Lemma \ref{lemma:joint.dens}, and the last equality follows from Lemma \ref{lemma:desint}. Hence, $\mu_{\SD(U)}\ll\lambda$ and the \RN derivative has the desired expression.
\end{myproof}

\begin{myproof}[Proof of Theorem \ref{theo:cond.dens}.] Let $\mu_{U,\SD(U)}$ be the joint distribution of $U$ and $\SD(U)$, and let $\mu_{\SD(U)}$ be the distribution of $\SD(U)$. For $A\in \cR^k$ and $B\in\cR$, by the definition of conditional distribution, we have that
\begin{align*}
\mu_{U,\SD(U)}(A\times B) &= \PP\{U\in A, \SD(U)\in B\} = \int_B \mu_{U\mid\SD(U)}(A\mid r)\,d\mu_{\SD(U)}(r) \\
                          &= \int_B \mu_{U\mid\SD(U)}(A\mid r)\,\frac{d\mu_{\SD(U)}}{d\lambda}(r)\,d\lambda(r) \, ,
\end{align*}
in which we have used the Leibniz rule for the \RN derivatives. On the other hand, by Lemmas \ref{lemma:joint.dens} and \ref{lemma:desint}, we have that
\begin{align*}
\mu_{U,\SD(U)}(A\times B) &= \int_{A\times B} f_U(u)\,\I_{\HH_r}(u)\,d\xi(u,r) \\
                          &= \int_{B}\left(\int_{A} f_U(u)\,\I_{\HH_r}(u)\,d\tau_r(u)\right)d\lambda(r) \, ,
\end{align*}
with $u\in\RR^k$ and $r\in\RR$. Both expressions for $\mu_{U,\SD(U)}(A\times B)$ are compatible if
\[
  \mu_{U\mid\SD(U)}(A\mid r) = \frac{{\displaystyle \int_A} f_U(u)\,\I_{\HH_r}(u)\,d\tau_r(u)}{f_{\SD(U)}(r)} \, ,
\]
for almost every $r$ $[\lambda]$. Therefore, we have that $\mu_{U\mid \SD(U)}(\;\cdot\mid r) \ll \tau_r$, for almost every $r>0$ $[\lambda]$, with \RN derivative $d\mu_{U\mid \SD(U)}/d\tau_r=f_{U\mid \SD(U)}(\:\cdot\mid r)$ given by
\[
  f_{U\mid \SD(U)}(u\mid r) = \frac{f_U(u)}{f_{\SD(U)}(r)}\,\I_{\HH_r}(u) \, ,
\]
as desired. The fact that $\mu_{U\mid\SD(U)}(\HH_r\mid r)=1$ follows immediately.
\end{myproof}

\begin{myproof}[Proof of Lemma \ref{lemma:model}.]
Let $\alpha_h$ be the measures over $(\RR^n,\cR^n)$ defined by $\alpha_h(A)=\int_A \left(\prod_{i=1}^k h_i^{c_i}\right)d\lambda_n(x)$, for each $h\in\HH_1$. Let $B=B_1\times\dots\times B_n$, with $B_i\in\cR$, for $i=1,\dots,n$. By the hypothesis of conditional independence and Tonelli's Theorem, we have that
\begin{align*}
  \mu_{X\mid H}(B\mid h) &= \prod_{j=1}^n \mu_{X_j\mid H}(B_j\mid h) = \prod_{j=1}^n \int_{B_j} f(x_j)\,d\lambda(x_j) \\
  &= \int_B \left( \prod_{j=1}^n f(x_j) \right)d\lambda_n(x) = \int_B \left( \prod_{j=1}^n \sum_{i=1}^k h_i\,\I_{[t_{i-1},t_i)}(x_j) \right) d\lambda_n(x) \\ &= \int_B \left( \prod_{i=1}^k h_i^{c_i} \right) d\lambda_n(x) =\alpha_h(B) \, .
\end{align*}
Hence, $\mu_{X\mid H}(\;\cdot\mid h)$ and $\alpha_h$ agree on the $\pi$-system of product sets that generate $\cR^n$. Therefore, by Theorem A.26 of \cite{schervish}, both measures agree on the whole sigma-field $\cR^n$. It follows that $\mu_{X\mid H}(\;\cdot\mid h)\ll\lambda_n$, almost surely $[\mu_H]$, and the \RN derivative has the desired expression.
\end{myproof}

\begin{myproof}[Proof of Theorem \ref{theo:closure}.]
By Bayes Theorem, for each $A\in\cR^k$, we have that
\begin{equation*}
\begin{split}
  \mu_{H\mid X}(A\mid x) &= C_0 \int_A f_{X\mid H}(x\mid h)\,d\mu_H(h) = C_0 \int_A \left(\prod_{i=1}^k h_i^{c_i}\right) d\mu_H(h) \\
  &= C_0 \int_A \left(\prod_{i=1}^k h_i^{c_i}\right) \frac{d\mu_H}{d\tau_1}(h)\,d\tau_1(h) \\
  &= \frac{C_0}{f_{\SD(U)}(1)} \int_A \left(\prod_{i=1}^k h_i^{c_i}\right) f_U(h)\,\I_{\HH_1}(h)\,d\tau_1(h) \, ,
\end{split}
\end{equation*}
in which we have used the expression of the likelihood obtained in Lemma \ref{lemma:model}, the Leibniz rule for the \RN derivatives, the expression of $d\mu_H/d\tau_1$ in Definition \ref{def:srd}, and the constant $C_0$ is such that $\mu_{H\mid X}(\HH_1\mid x)=1$. The remainder of the proof relies on some matrix algebra. Let $I$ be the identity matrix. Since,  by definition, $\Sigma$ is symmetric , we have that $I = I^\top = (\Sigma \Sigma^{-1})^\top = (\Sigma^{-1})^\top\Sigma^\top = (\Sigma^{-1})^\top\Sigma$. Therefore, we have that $(\Sigma^{-1})^\top=\Sigma^{-1}$. Write $l=\log h$. Since the scalar $l^\top\Sigma^{-1}m$ is equal to its transpose $(l^\top\Sigma^{-1}m)^\top=m^\top\Sigma^{-1}l$, we have that 
$$
  (l-m)^\top\Sigma^{-1}(l-m)=l^\top\Sigma^{-1}l -2m^\top\Sigma^{-1}l+m^\top\Sigma^{-1}m \, .
$$
Defining $d=\Sigma\,c$, we have
\begin{align*}
  & \left(\prod_{i=1}^k  h_i^{c_i}\right) \, \exp\left( -\frac{1}{2}(l-m^*)^\top\Sigma^{-1}(l-m^*)\right) \\
  & \qquad = \exp\left( -\frac{1}{2}\left( -2d^\top\Sigma^{-1}l +l^\top\Sigma^{-1}l -2m^\top\Sigma^{-1}l+m^\top\Sigma^{-1}m\right)\right) \\
  & \qquad = C_1 \, \exp\biggl( -\frac{1}{2}\left( -2d^\top\Sigma^{-1}l +l^\top\Sigma^{-1}l-2m^\top\Sigma^{-1}l+m^\top\Sigma^{-1}m\right) \\ & \qquad\qquad\qquad\qquad\qquad\qquad\qquad\qquad\qquad\qquad +2m^\top\Sigma^{-1}d + d^\top\Sigma^{-1}d\biggr) \, ,
\end{align*}
with $C_1=\exp\left(-(1/2)\left(-2m^\top\Sigma^{-1}d-d^\top\Sigma^{-1}d\right)\right)$. Define $m^*=m+d$. Since the scalar $d^\top\Sigma^{-1}m=(d^\top\Sigma^{-1}m)^\top=m^\top\Sigma^{-1}d$, we have that $(m^*)^\top\Sigma^{-1}m^*=m^\top\Sigma^{-1}m +2m^\top\Sigma^{-1}d+d^\top\Sigma^{-1}d$. Hence, we obtain
\begin{align*}
  \left(\prod_{i=1}^k h_i^{c_i}\right) \, &\exp\left( -\frac{1}{2}(l-m^*)^\top\Sigma^{-1}(l-m^*)\right) \\
  &= C_1 \, \exp\left( -\frac{1}{2}\left( l^\top\Sigma^{-1}l -2(m^*)^\top\Sigma^{-1}l+(m^*)^\top\Sigma^{-1}m^*\right)\right) \\
  &= C_1 \, \, \exp\left( -\frac{1}{2}(l-m^*)^\top\Sigma^{-1}(l-m^*)\right) \, .
\end{align*}
Using this result in the expression of $\mu_{H\mid X}$ together with the expression of $f_U$, we have
\[
  \mu_{H\mid X}(A\mid x) = C_2 \int_A f_{U^*}(h)\,\I_{\HH_1}(h)\,d\tau_1(h) \, ,
\]
in which $C_2=(C_0\,C_1)/f_{\SD(U)}(1)$ and $f_{U^*}$ is a density of the random vector $U^*\sim L_k(m^*,\Sigma)$. We conclude that, given that $X=x$, the vector $H$ has the distribution of the heights of the steps of $\phi^*\sim\Delta(m^*,\Sigma)$, as desired.
\end{myproof}

\begin{proposition}\label{prop:pred}
Suppose that the random variables $X_1,\dots,X_{n+1}$ are modeled according to Lemma \ref{lemma:model}. Denote by $\mu_{X_i}$ the distribution of $X_i$, for $i=1,\dots,n+1$. For convenience, use the notations $X^{(n)}=(X_1,\dots,X_n)$ and $x^{(n)}=(x_1,\dots,x_n)\in\RR^k$. Then, for every $A\in\cR$, we have

(a) $\mu_{X_i}(A) = {\displaystyle\int_A} \E[\phi(y)]\,d\lambda(y)$, for $i=1,\dots,n+1$;

(b) $\mu_{X_{n+1}\mid X^{(n)}}(A\mid x^{(n)}) = {\displaystyle\int_A} \E[\phi(y)\mid X^{(n)}=x^{(n)}]\,d\lambda(y)$, a.s. $[\mu_{X^{(n)}}]$.
\end{proposition}

\begin{myproof}
By Definition \ref{def:srd}, we have
\[
  \E[\phi(y)] = \E\left[ \sum_{i=1}^k H_i\,\I_{[t_{i-1},t_i)}(y) \right] = \int_{\RR^k} f(y)\,d\mu_H(h) \, ,
\]
in which $h\in\RR^k$ and $f(y) = \sum_{i=1}^k h_i\,\I_{[t_{i-1},t_i)}(y)$, for $y\in\RR$. In an analogous manner, we have
\[
  \E[\phi(y)\mid X^{(n)}=x^{(n)}] = \int_{\RR^k} f(y)\,d\mu_{H\mid X^{(n)}}(h\mid x^{(n)}) \, .
\]
For item (a), note that
\begin{align*}
  \mu_{X_i}(A) &= \PP\{X_i\in A, H\in\RR^k\} = \int_{\RR^k} \mu_{X_i\mid H}(A\mid h)\,d\mu_H(h) \\
  &= \int_{\RR^k} \left(\int_A f(y)\,d\lambda(y)\right) d\mu_H(h) = \int_A \left(\int_{\RR^k} f(y)\,d\mu_H(h)\right) d\lambda(y) \\ 
  &= \int_A \E[\phi(y)]\,d\lambda(y) \, ,
\end{align*}
in which the fourth equality follows from Tonelli's Theorem. For item (b), for each $B\in\cR^n$, we have
\[
  \PP\{X_{n+1}\in A, X^{(n)}\in B\} = \int_B \mu_{X_{n+1}\mid X^{(n)}}(A\mid x^{(n)})\,d\mu_{X^{(n)}}(x^{(n)}) \, .
\]
On the other hand, we have
\begin{align*}
  \PP\{X_{n+1}\in A, &X^{(n)}\in B\} = \PP\{X_{n+1}\in A, X^{(n)}\in B, H\in\RR^k\} \\
  &= \int_{B\times\RR^k} \mu_{X_{n+1}\mid X^{(n)},H}(A\mid x^{(n)},h)\,d\mu_{X^{(n)},H}(x^{(n)},h) \\
  &= \int_{B\times\RR^k} \mu_{X_{n+1}\mid H}(A\mid h)\,d\mu_{X^{(n)},H}(x^{(n)},h) \\
  &= \int_B \left( \int_{\RR^k} \mu_{X_{n+1}\mid H}(A\mid h)\,d\mu_{H\mid X^{(n)}}(h\mid x^{(n)}) \right) d\mu_{X^{(n)}}(x^{(n)}) \\
  &= \int_B \left( \int_{\RR^k} \left( \int_A f(y)\,d\lambda(y) \right) d\mu_{H\mid X^{(n)}}(h\mid x^{(n)}) \right) d\mu_{X^{(n)}}(x^{(n)}) \\
  &= \int_B \left( \int_A \left( \int_{\RR^k} f(y)\,d\mu_{H\mid X^{(n)}}(h\mid x^{(n)}) \right) d\lambda(y) \right) d\mu_{X^{(n)}}(x^{(n)}) \\
  &= \int_B \left( \int_A \E[\phi(y)\mid X^{(n)}=x^{(n)}]\,d\lambda(y) \right) d\mu_{X^{(n)}}(x^{(n)}) \, ,
\end{align*}
in which the third equality follows from the hypothesis of conditional independence and Theorem B.61 of \cite{schervish}, the fourth equality is a consequence of Theorem 2.6.4 of \cite{ash}, and the sixth equality is due to Tonelli's Theorem. Comparing both expressions for $\PP\{X_{n+1}\in A, X^{(n)}\in B\}$, we get the desired result.
\end{myproof}

\begin{proposition}\label{prop:Lkr}
Let $\mu_K=P\circ K^{-1}$ over $(\NN,2^\NN)$ be the distribution of $K$ and let $\mu_R=P\circ R^{-1}$ over $(\RR,\cR)$be the distribution of $R$. Denote by $\mu_{K,R}$ the joint distribution of $K$ and $R$, which by the independence of $K$ and $R$ is equal to the product measure $\mu_K\times\mu_R$, and let $\mu_{K,R,H}$ be the joint distribution of $K$, $R$ and $H$. In the hierarchical model described in Section \ref{sec:rp}, we have that $\mu_{X\mid K,R}(\;\cdot\mid k,\rho)\ll\lambda_n$, almost surely $[\mu_{K,R}]$, with \RN derivative
\[
  \frac{d\mu_{X\mid K,R}}{d\lambda_n}(x\mid k,\rho) = f_{X\mid K,R}(x\mid k,\rho) =  \int_{\RR^k} f_{X\mid H}(x\mid h)\,d\mu_{H\mid K,R}(h\mid k,\rho) \, ,
\]
for the $f_{X\mid H}$ defined on Lemma \ref{lemma:model}.
\end{proposition}

\begin{myproof}
Let $A\in\cR^n$ and $B\in 2^\NN\otimes\cR$. By the definition of conditional distribution, we have
\[
  \PP\{X\in A, (K,R)\in B\} = \int_B \mu_{X\mid K,R}(A\mid k,\rho)\,d\mu_{K,R}(k,\rho) \, .
\]
On the other hand, by arguments similar to those used in the proof of Proposition \ref{prop:pred}, we have
\begin{align*}
  \PP\{X\in &A, (K,R)\in B\} = \PP\{X\in A, (K,R)\in B, H\in\RR^k\} \\
  &= \int_{B\times\RR^k} \mu_{X\mid K,R,H}(A\mid k,\rho,h)\,d\mu_{K,R,H}(k,\rho,h) \\
  &= \int_{B\times\RR^k} \mu_{X\mid H}(A\mid h)\,d\mu_{K,R,H}(k,\rho,h) \\
  &= \int_B \left( \int_{\RR^k} \mu_{X\mid H}(A\mid h)\,d\mu_{H\mid K,R}(h\mid k,\rho) \right) d\mu_{K,R}(k,\rho) \\
  &= \int_B \left( \int_{\RR^k} \left( \int_A f_{X\mid H}(x\mid h)\,d\lambda_n(x) \right) d\mu_{H\mid K,R}(h\mid k,\rho) \right) d\mu_{K,R}(k,\rho) \\
  &= \int_B \left( \int_A \left( \int_{\RR^k} f_{X\mid H}(x\mid h)\,d\mu_{H\mid K,R}(h\mid k,\rho) \right) d\lambda_n(x) \right) d\mu_{K,R}(k,\rho) \, .
\end{align*}
Comparing both expressions for $\PP\{X\in A, (K,R)\in B\}$, we have
\[
\mu_{X\mid K,R}(A\mid k,\rho) = \int_A \left( \int_{\RR^k} f_{X\mid H}(x\mid h)\,d\mu_{H\mid K,R}(h\mid k,\rho) \right) d\lambda_n(x) \, ,
\]
almost surely $[\mu_{K,R}]$, and the result follows.
\end{myproof}

\begin{myproof}[Proof of Proposition \ref{prop:decision}.]
By Tonelli's Theorem, the expected loss is
\[
  \E[L(\phi,f)] = \int_a^b f^2(x)\,d\lambda(x) - 2 \int_a^b f(x) \E[\phi(x)]\,d\lambda(x) + C_0 \, ,
\]
in which we have defined the positive constant $C_0=\int_a^b \E[\phi^2(x)]\,d\lambda(x)$. By hypothesis, each $f$ has the form $f(x) = \sum_{i=1}^N \alpha_i\,g_i(x)$, leading us to 
\begin{align*}
  \E[L(\phi,f)] &= \sum_{i,j=1}^N \left( \alpha_i \alpha_j \int_a^b g_i(x) g_j(x)\,d\lambda(x) \right) \\
  &- 2\sum_{i=1}^N \left( \alpha_i \int_a^b g_i(x) \E[\phi(x)]\,d\lambda(x) \right) + C_0 \, ,
\end{align*}
in which we have used the linearity of the integral. Therefore, minimizing the expected loss is the same as solving the problem of constrained minimization of the quadratic form $Q$. For the matrix $M=(M_{ij})$, note that, for every non null $y=(y_1,\dots,y_N)^\top\in\RR^N$, we have
\begin{align*}
  y^\top M y &= \sum_{i,j=1}^N y_i y_j M_{ij} = \sum_{i,j=1}^N \left( y_i y_j \int_a^b g_i(x) g_j(x)\,d\lambda(x) \right) \\
  &= \int_a^b \sum_{i,j=1}^N \left( y_i\,g_i(x)\,y_j\,g_j(x) \right) \,d\lambda(x) = \int_a^b \left( \sum_{i=1}^N  y_i\,g_i(x) \right)^2 \,d\lambda(x) > 0 \, , 
\end{align*}
in which we have used the linearity of the integral. Therefore, the matrix $M$ is positive definite, yielding (see \cite{bazaraa}) that the quadratic form $Q$ is convex and the problem of constrained minimization of $Q$ has a single global solution $(\hat{\alpha_1},\dots,\hat{\alpha_N})$. Since the Bayes decision is the $f$ that minimizes the expected loss, the result follows.
\end{myproof}

\bibliographystyle{plain}
\bibliography{srd}

\end{document}